\documentclass[12 pt,twoside]{amsart}

\usepackage{amsopn}
\usepackage{amssymb}
\usepackage{amscd}

\newtheorem{theorem}{Theorem}[section]
\newtheorem{lemma}[theorem]{Lemma}

\newtheorem{corollary}[theorem]{Corollary}

\newtheorem*{th1}{Theorem}

\theoremstyle{definition}
\newtheorem{definition}[theorem]{Definition}

\theoremstyle{remark}
\newtheorem{remark}[theorem]{Remark}

\numberwithin{equation}{section}

\begin{document}

\title[$J$-class weighted shifts]{$J$-class weighted shifts on the space of bounded sequences of complex numbers}

\author[George Costakis]{George Costakis}
\address{Department of Mathematics, University of Crete, Knossos Avenue, GR-714 09 Heraklion, Crete, Greece}
\email{costakis@math.uoc.gr}
\thanks{}

\author[Antonios Manoussos]{Antonios Manoussos}
\address{Fakult\"{a}t f\"{u}r Mathematik, SFB 701, Universit\"{a}t Bielefeld, Postfach 100131, D-33501 Bielefeld, Germany}
\email{amanouss@math.uni-bielefeld.de}
\thanks{During this research the second author was fully supported by SFB 701 ``Spektrale Strukturen und
Topologische Methoden in der Mathematik" at the University of Bielefeld, Germany. He would also like to express
his gratitude to Professor H. Abels for his support.}

\subjclass[2000]{Primary 47A16; Secondary 37B99, 54H20}

\date{}

\keywords{Hypercyclic operators, $J$-class operators, $J^{mix}$-class operators, unilateral and bilateral
weighted shifts.}

\begin{abstract}
We provide a characterization of $J$-class and $J^{mix}$-class unilateral weighted shifts on
$l^{\infty}(\mathbb{N})$ in terms of their weight sequences. In contrast to the previously mentioned result we
show that a bilateral weighted shift on $l^{\infty}(\mathbb{Z})$ cannot be a $J$-class operator.
\end{abstract}

\maketitle

\section{Introduction}
During the last years the dynamics of linear operators on infinite dimensional spaces has been extensively
studied, see the survey articles \cite{BoMaPe}, \cite{GE}, \cite{GE2}, \cite{GE3}, \cite{MoSa2}, \cite{Sh} and
the recent book \cite{BM}. Let us recall the notion of \textit{hypercyclicity}. Let $X$ be a separable Banach
space and $T:X\to X$ be a bounded linear operator. The operator $T$ is said to be hypercyclic provided there
exists a vector $x\in X$ such that its orbit under $T$, $Orb(T,x)=\{ T^nx : n=0,1,2,\ldots \}$, is dense in $X$.
If $X$ is Banach space (possibly non-separable) and $T:X\to X$ is a bounded linear operator then $T$ is called
\textit{topologically transitive} (\textit{topologically mixing}) if for every pair of non-empty open subsets
$U, V$ of $X$ there exists a positive integer $n$ such that $T^nU\cap V\neq\emptyset$ ($T^mU\cap V\neq\emptyset$
for every $m\geq n$ respectively). It is well known, and easy to prove, that if $T$ is a bounded linear operator
acting on separable Banach space $X$ then $T$ is hypercyclic if and only if $T$ is topologically transitive.

A first step to understand the dynamics of linear operators is to look at particular operators as for example
the weighted shifts. Salas \cite{Salas2} was the first who characterized the hypercyclic weighted shifts in
terms of their weight sequences. We would like to point out that $l^{\infty}(\mathbb{N})$ and
$l^{\infty}(\mathbb{Z})$ do not support hypercyclic operators since they are not separable Banach spaces. In
fact they do not support topologically transitive operators as it was shown by Berm\'{u}dez and Kalton in
\cite{BeKa}. Recently B\`{e}s, Chan and Sanders \cite{bcs} showed that there exists a weak* hypercyclic weighted
shift $T$ on $l^{\infty}(\mathbb{N})$, i.e there exists a vector $x\in l^{\infty}(\mathbb{N})$ whose orbit
$Orb(T,x)$ is dense in the weak* topology of $l^{\infty}(\mathbb{N})$. In fact they give a characterization of
the weak* hypercyclic weighted shifts in terms of their weight sequences. In \cite{cm} we studied the dynamics
of operators by replacing the orbit of a vector with its extended limit set. To be precise, let $T:X\to X$ be a
bounded linear operator on a Banach space $X$ (not necessarily separable) and $x\in X$. A vector $y$ belongs to
the extended limit set $J(x)$ of $x$ if there exist a strictly increasing sequence of positive integers $\{
k_n\}$ and a sequence $\{ x_n\}\subset X$ such that $x_n\to x$ and $T^{k_n}x_n\to y$. If $J(x)=X$ for some
non-zero vector $x\in X$ then $T$ is called $J$-class operator. Roughly speaking, the use of the extended limit
set ``localizes" the notion of hypercyclicity. The last can be justified by the following: $J(x)=X$ if and only
if for every open neighborhood $U$ of $x$ and every non-empty open set $V\subset X$ there exists a positive
integer $n$ such that $T^nU\cap V\neq\emptyset$.

The purpose of this paper is to study the dynamical behavior of weighted shifts on the spaces of bounded
sequences of complex numbers $l^{\infty}(\mathbb{N})$ and $l^{\infty}(\mathbb{Z})$ through the use of the
extended limit sets. Our main result is the following (see Theorem \ref{t1}).

\begin{th1}
Let $T:l^{\infty}(\mathbb{N})\to l^{\infty}(\mathbb{N})$ be a backward unilateral weighted shift with positive
weights $(\alpha_n )_{n\in\mathbb{N}}$. The following are equivalent.
\begin{enumerate}
\item[(i)] $T$ is a $J$-class operator.

\item[(ii)] $\displaystyle{\lim_{n\to +\infty} \left( \inf_{j\geq 0}\prod_{i=1}^{n} \alpha_{i+j}  \right)
=+\infty.}$
\end{enumerate}
In particular, if $T$ is a $J$-class operator then the sequence of weights $(\alpha_n )_{n\in\mathbb{N}}$ is bounded from below by a positive number and we have the
following complete description of the set of $J$-vectors.
\[
\{ x\in l^{\infty}(\mathbb{N}): J(x)=l^{\infty}(\mathbb{N}) \}=c_0(\mathbb{N}),
\]
where $c_0(\mathbb{N})=\{ x=(x_n)_{n\in \mathbb{N}}\in l^{\infty}(\mathbb{N}): \lim_{n\to +\infty}x_n= 0 \}$.
\end{th1}

Observe that if $T$ is a $J$-class backward unilateral weighted shift on $l^{\infty}(\mathbb{N})$ then in view
of the above theorem and Salas' characterization of hypercyclic weighted shifts, see \cite{Salas2}, we conclude
that $T$ is hypercyclic on $l^{p}(\mathbb{N})$ for every $1\leq p<+\infty$. However, as we show in section 3,
the converse is not always true.

On the other hand the situation is completely different in the case of bilateral weighted shifts. In particular
we show that a bilateral weighted shift on $l^{\infty}(\mathbb{Z})$ cannot be a $J$-class operator, see Theorem
\ref{t2}. In addition, we prove similar results for $J^{mix}$-class weighted shifts (see Definitions 2.1 and
2.2).

\section{Preliminaries}
\begin{definition}
Let $T:X\to X$ be a bounded linear operator on a Banach space $X$. For every $x\in X$ the sets
\[
\begin{split}
J(x)=\{ &y\in X:\,\mbox{ there exist a strictly increasing sequence of positive}\\
&\mbox{integers}\,\{k_{n}\}\,\mbox{and a sequence }\,\{x_{n}\}\subset X\,\mbox{such that}\, x_{n}\rightarrow x\,
\mbox{and}\\
&T^{k_{n}}x_{n}\rightarrow y\},
\end{split}
\]
\[
\begin{split}
J^{mix}(x)=\{ &y\in X:\,\mbox{ there exists a sequence}\,\{x_{n}\}\subset X\,\mbox{such that}\\
&x_{n}\rightarrow x\,\,\mbox{and}\,\,T^{n}x_{n}\rightarrow y\}
\end{split}
\]
will be called the extended limit set of $x$ under $T$ and the extended mixing limit set of $x$ under $T$
respectively.
\end{definition}

\begin{definition}
A bounded linear operator $T:X\to X$ acting on a Banach space $X$ will be called a $J$-class ($J^{mix}$-class)
operator if there exists a non-zero vector $x\in X$ such that $J(x)=X$ ($J^{mix}(x)=X$ respectively).
\end{definition}

\begin{definition}
Let $T$ be a bounded linear operator acting on a Banach space $X$. A vector $x\in X$ will be called a $J$-vector
($J^{mix}$-vector) if $J(x)=X$ ($J^{mix}(x)=X$ respectively).
\end{definition}

\begin{remark}
Observe that
\begin{enumerate}
\item[(i)] an operator $T:X\to X$ is topologically transitive if and only if $J(x)=X$ for every $x\in X$,

\item[(ii)] an operator $T:X\to X$ is topologically mixing if and only if $J^{mix}(x)=X$ for every $x\in X$,
\end{enumerate}
see \cite{cm}. Hence every hypercyclic operator (topologically mixing) is a $J$-class operator ($J^{mix}$-class
operator). However the converse is not true. To see that consider the operator $3I\oplus 2B: \mathbb{C}\oplus
l^{2}(\mathbb{N})\to \mathbb{C}\oplus l^{2}(\mathbb{N})$ where $I$ is the identity map on $\mathbb{C}$ and $B$
is the backward shift on the space of square summable sequences $l^{2}(\mathbb{N})$. Consider any non-zero
vector $x\in l^{2}(\mathbb{N})$.  We shall prove that $J_{3I\oplus 2B}^{mix}(0\oplus x)=\mathbb{C}\oplus
l^{2}(\mathbb{N})$. Let $y\in l^{2}(\mathbb{N})$ and $\lambda\in\mathbb{C}$. There exists a sequence  $\{ x_n\}$
in $l^{2}(\mathbb{N})$ such that $(2B)^{n}x_n\to y$. Define the vectors $\frac{\lambda}{3^{n}}\oplus x_n$. Then
we have $\frac{\lambda}{3^{n}}\oplus x_n\to 0\oplus x$ and $(3I\oplus 2B)^{n}(\frac{\lambda}{3^{n}}\oplus
x_n)\to \lambda\oplus y$. Hence $3I\oplus 2B$ is a $J^{mix}$-class operator which is not hypercyclic. In fact it
is not even supercyclic, see \cite{GLSR}.

Let us also give an example of a backward weighted shift, acting on a non-separable space, which is a $J$-class
operator but not topologically transitive. Consider the operator $2B:l^{\infty}(\mathbb{N})\to
l^{\infty}(\mathbb{N})$ where $B$ is the backward shift and $l^{\infty}(\mathbb{N})$ is the space of bounded
sequences. Theorem \ref{t3} implies that $2B$ is a $J^{mix}$-class operator. On the other hand the space
$l^{\infty}(\mathbb{N})$ does not support topologically transitive operators, see \cite{BeKa}.
\end{remark}

The next lemma, which will be of use to us, also appears in \cite{cm}. For the convenience of the reader we give
its proof.

\begin{lemma}\label{l1}
Let $T:X\to X$ be a bounded linear operator on a Banach space $X$ and $\{ x_{n}\}$, $\{ y_{n} \}$ be two
sequences in $X$ such that $x_{n}\rightarrow x$ and $y_{n} \rightarrow y$ for some $x,y\in X$.
\begin{enumerate}
\item[(i)] If $y_{n} \in J(x_{n})$ for every $n=1,2,\ldots$, then $y\in J(x)$.

\item[(ii)] If $y_{n} \in J^{mix}(x_{n})$ for every $n=1,2,\ldots$, then $y\in J^{mix}(x)$.
\end{enumerate}
\end{lemma}
\begin{proof}
(i) For $n=1$ there exists a positive integer $k_{1}$ such that
$$
\Vert x_{k_{1}}-x\Vert<\frac{1}{2}\,\,\mbox{and}\,\,\Vert y_{k_{1}} -y\Vert<\frac{1}{2}.
$$
Since $y_{k_{1}}\in J(x_{k_{1}})$ we may find a positive integer $l_{1}$ and $z_{1}\in X$ such that
$$
\Vert z_{1}-x_{k_{1}}\Vert<\frac{1}{2}\,\,\mbox{and}\,\,\Vert T^{l_{1}} z_{1}-y_{k_{1}}\Vert<\frac{1}{2}.
$$
Therefore,
$$
\Vert z_{1}-x\Vert<1\,\,\mbox{and}\,\,\Vert T^{l_{1}}z_{1}-y\Vert<1.
$$
Proceeding inductively we find a strictly increasing sequence of positive integers $\{l_{n}\}$ and a sequence
$\{z_{n}\}$ in $X$ such that
$$
\Vert z_{n}-x\Vert<\frac{1}{n}\,\,\mbox{and}\,\,\Vert T^{l_{n}}z_{n} -y\Vert<\frac{1}{n}.
$$
This completes the proof of assertion (i).

(ii) For $n=1$ there exists a positive integer $k_{1}$ such that
$$
\Vert x_{k_{1}}-x\Vert<\frac{1}{2}\,\,\mbox{and}\,\,\Vert y_{k_{1}} -y\Vert<\frac{1}{2}.
$$
Since $y_{k_{1}}\in J^{mix}(x_{k_{1}})$ we may find a positive integer $l_{1}$ and a sequence $\{ z_{n}\}\subset
X$ such that
$$
\Vert z_{n}-x_{k_{1}}\Vert<\frac{1}{2}\,\,\mbox{and}\,\,\Vert T^{n} z_{n}-y_{k_{1}}\Vert<\frac{1}{2}
$$
for every $n\geq l_{1}$. Therefore,
$$
\Vert z_{n}-x\Vert<1\,\,\mbox{and}\,\,\Vert T^{n}z_{n}-y\Vert<1
$$
for every $n\geq l_{1}$. Proceeding  in the same way we may find a positive integer $l_{2}>l_{1}$ and a sequence
$\{ w_{n}\}\subset X$ such that
$$
\Vert w_{n}-x\Vert<\frac{1}{2}\,\,\mbox{and}\,\,\Vert T^{n}w_{n}-y\Vert<\frac{1}{2}
$$
for every $n\geq l_{2}$. Set $v_{n}=z_{n}$ for every $l_{1}\leq n<l_{2}$, hence
$$
\Vert v_{n}-x\Vert<1\,\,\mbox{and}\,\,\Vert T^{n}v_{n}-y\Vert<1.
$$
Proceeding inductively we find a strictly increasing sequence of positive integers $\{ n_k\}$ and a sequence $\{v_{n}\}$ in $X$ such that if $n\geq n_k$ then
$$
\Vert v_{n}-x\Vert<\frac{1}{k}\,\,\mbox{and}\,\,\Vert T^{n}v_{n} -y\Vert<\frac{1}{k}.
$$
Take any $\epsilon >0$. There exists a positive integer $k_0$ such that $\frac{1}{k_0}<\epsilon$. Hence for every $n\geq n_{k_0}$ we get
$$
\Vert v_{n}-x\Vert<\frac{1}{k_0}<\epsilon\,\,\mbox{and}\,\,\Vert T^{n}v_{n} -y\Vert<\frac{1}{k_0}<\epsilon.
$$
This completes the proof of assertion (ii).
\end{proof}

\section{Main results}
\begin{theorem}\label{t1}
Let $T:l^{\infty}(\mathbb{N})\to l^{\infty}(\mathbb{N})$ be a backward unilateral weighted shift with positive
weights $(\alpha_n )_{n\in\mathbb{N}}$. The following are equivalent.
\begin{enumerate}
\item[(i)] $T$ is a $J$-class operator.

\item[(ii)] $\displaystyle{\lim_{n\to +\infty} \left( \inf_{j\geq 0}\prod_{i=1}^{n} \alpha_{i+j}  \right) =+\infty.}$
\end{enumerate}
In particular, if $T$ is a $J$-class operator then the sequence of weights $(\alpha_n )_{n\in\mathbb{N}}$ is bounded from below by a positive number and we have the
following complete description of the set of $J$-vectors.
\[
\{ x\in l^{\infty}(\mathbb{N}): J(x)=l^{\infty}(\mathbb{N}) \}=c_0(\mathbb{N}),
\]
where $c_0(\mathbb{N})=\{ x=(x_n)_{n\in \mathbb{N}}\in l^{\infty}(\mathbb{N}): \lim_{n\to +\infty}x_n= 0 \}$.
\end{theorem}
\begin{proof}
Let us prove that (i) implies (ii). There exists a non-zero vector $x\in l^{\infty}(\mathbb{N})$ such that
$J(x)=l^{\infty}(\mathbb{N})$. Consider the vector $y=(1,1,\ldots )$. Then there exists a strictly increasing
sequence $\{ k_n\}$ of positive integers and a sequence $\{ y_n\}\in l^{\infty}(\mathbb{N})$,
$y_n=(y_{nm})_{m=1}^{\infty}$,  such that
\[
\| y_n -x\|_{\infty} \to 0 \quad \textrm{and}\quad \| T^{k_n}y_n - (1,1,\ldots )\|_{\infty} \to 0.
\]
Observe that
\[
\| T^{k_n}y_{n} - (1,1,\ldots )\|_{\infty} = \sup_{j\geq 0} \left| \left (\prod_{i=1}^{k_n} \alpha_{i+j}\right ) y_{n(k_n+j+1)} -1 \right| \to 0
\]
as $n\to\infty$. Fix $0<\epsilon <1$. There exists a positive integer $n_1$ such that
\begin{equation}\label{e1}
\| y_n-x\|_{\infty} <\epsilon \quad\textrm{for every}\quad n\geq n_1
\end{equation}
and
\[
\sup_{j\geq 0} \left| \left ( \prod_{i=1}^{k_{n_1}}\alpha_{i+j}\right ) y_{n_1(k_{n_1}+j+1)} -1 \right| <\epsilon .
\]
Therefore
\begin{equation}\label{e2}
\left| \left ( \prod_{i=1}^{k_{n_1}} \alpha_{i+j} \right ) y_{n_1(k_{n_1}+j+1)} \right| >1-\epsilon \quad\textrm{for every}\quad j\geq 0.
\end{equation}
On the other hand, using (\ref{e1}), we have
\begin{equation}\label{e3}
\begin{split}
\left| \left ( \prod_{i=1}^{k_{n_1}} \alpha_{i+j}\right ) y_{n_1(k_{n_1}+j+1)} \right| &\leq  \left ( \prod_{i=1}^{k_{n_1}}
\alpha_{i+j} \right ) \|y_{n_1}\|_{\infty}\\
&< \left ( \prod_{i=1}^{k_{n_1}} \alpha_{i+j}  \right ) (\epsilon +\| x\|_{\infty})
\end{split}
\end{equation}
for every $j\geq 0$. By (\ref{e2}) and (\ref{e3}) it follows that
\[
\prod_{i=1}^{m_1} \alpha_{i+j}  > \frac{1-\epsilon}{\epsilon +\| x\|_{\infty}}\quad\textrm{for every}\quad j\geq 0,
\]
where $m_1:= k_{n_1}$. For every $l=2,3,\ldots$ consider the vector $(l,l,\ldots )$. Since $(l,l,\ldots )\in
J(x)$ and working as before we inductively construct a strictly increasing sequence $\{ m_l\}$ of positive
integers such that
\[
\prod_{i=1}^{m_l} \alpha_{i+j} > \frac{l-\epsilon}{\epsilon +\| x\|_{\infty}} \quad\textrm{for every}\quad j\geq 0\quad\textrm{and every}\quad l\geq 1.
\]
The last implies that
\[
\lim_{l\to +\infty} \left( \inf_{j\geq 0}\prod_{i=1}^{m_l} \alpha_{i+j}  \right) =+\infty
\]
which in turn yields
\[
\limsup_{n\to +\infty} \left( \inf_{j\geq 0}\prod_{i=1}^{n} \alpha_{i+j}  \right) =+\infty.
\]
It remains to show that
\[
\lim_{n\to +\infty} \left( \inf_{j\geq 0}\prod_{i=1}^{n} \alpha_{i+j}  \right) =+\infty.
\]
Let us first show that the sequence $(\alpha_n )_{n\in\mathbb{N}}$ is bounded from below by a positive number. Fix a positive number $M>1$. There exists a positive
integer $N$ such that
\[
\prod_{i=1}^{N} \alpha_{i+j} >M \quad\mbox{for every}\,\, j\geq 0.
\]
If $N=1$ there is nothing to prove. Assume that $N>1$. For every $j\geq 0$ and since $\| T\| =\sup_n \alpha_n$, we have
\[
\alpha_{j+1} \| T\|^{N-1} \geq \alpha_{j+1} \left ( \prod_{i=2}^{N} \alpha_{i+j}\right )>M.
\]
Proceeding inductively we conclude that
\[
\alpha_n \geq \frac{M}{\| T\|^{N-1}}
\]
for every $n\in\mathbb{N}$. Take any positive integer $n>N$. There exist positive integers $p_n, v_n$ such that
$n=Np_n+v_n$ and $0\leq v_n\leq N-1$. Since $( \alpha_n)_{n\in\mathbb{N}}$ is bounded from below by $\frac{M}{\|
T\|^{N-1}}$ it follows that
\[
\prod_{i=1}^{n} \alpha_{i+j} > M^{p_n}\,C  \quad\mbox{for every}\,\, j\geq 0,
\]
where
\[
C=\min\left \{ \left ( \frac{M}{\| T\|^{N-1}} \right)^{N-1},1 \right\}.
\]
From the last and the fact that $M>1$ it clearly follows that
\[
\lim_{n\to +\infty} \left ( \inf_{j\geq 0} \prod_{i=1}^{n}\alpha_{i+j}\right ) =+\infty.
\]

We shall now prove that (ii) implies (i). Fix a vector $x=(x_1,x_2,\ldots )$ in $l^{\infty}(\mathbb{N})$ with finite support. There exists a positive integer $n_0$
such that $x_n=0$ for every $n\geq n_0$ and $\inf_{j\geq 0} \prod_{i=1}^{n}\alpha_{i+j}>0$ for every $n\geq n_0$. Consider any vector $y=(y_1,y_2,\ldots )\in
l^{\infty}(\mathbb{N})$. We set
\[
y_n=\left(x_1,x_2,\ldots ,x_{n_0-1},0,\ldots, 0,\frac{y_1}{\prod_{i=1}^{n} \alpha_{i}},\frac{y_2}{\prod_{i=1}^{n} \alpha_{i+1}},\frac{y_3}{\prod_{i=1}^{n}
\alpha_{i+2}},\ldots \right)
\]
for every $n\geq n_0$, where the $0$'s  fill all the coordinates from the $n_0$-th up to $n$-th position. Then for every $n\geq n_0$ we have
\[
\| y_n-x\|_{\infty}=\sup_{j\geq 0}\left|\frac{y_{j+1}}{\prod_{i=1}^{n} \alpha_{i+j}} \right| \leq \frac{\| y\|_{\infty}}{\inf_{j\geq 0} \prod_{i=1}^{n} \alpha_{i+j}},
\]
hence $y_n\to x$. Observe also that $T^{n}y_n=y$, so $y\in J(x)$. Thus $T$ is a $J$-class operator and this completes the proof that (ii) implies (i).

It remains to show that the set of $J$-vectors is $c_0(\mathbb{N})$. From the proof that (ii) implies (i) we
have that if $x$ is a vector with finite support then $J(x)=l^{\infty}(\mathbb{N})$. Since the closure of the
set of all vectors with finite support is $c_{0}(\mathbb{N})$, by Lemma \ref{l1}, we conclude that
\[
c_0(\mathbb{N})\subset\{ x\in l^{\infty}(\mathbb{N}): J(x)=l^{\infty}(\mathbb{N}) \}.
\]
To prove the converse inclusion, take a vector $x$ such that $J(x)=l^{\infty}(\mathbb{N})$. Consider the zero vector  and let $\epsilon$ be a positive number. There
exist positive integers $n_0,n_1$ and a vector $y_{n_0}=(y_{n_0k})_{k\in\mathbb{N}}$ such that
\[
\| y_{n_0}-x\|_{\infty} <\epsilon,\,\, \| T^{n_1}y_{n_0}\|_{\infty}<\epsilon \quad\mbox{and}\quad \prod_{i=1}^{n_1} \alpha_{i+j}>1\,\,\,\mbox{for every}\,\, j\geq 0.
\]
Hence we have
\[
\left| \left ( \prod_{i=1}^{n_1} \alpha_{i+j} \right ) y_{n_{0}(n_1+j+1)} \right| <\epsilon
\]
for every $j\geq 0$. The last and the previous bound on the weights imply that
\[
| y_{n_{0}(n_1+j+1)} |<\frac{\epsilon}{\prod_{i=1}^{n_1} \alpha_{i+j}}<\epsilon
\]
for every $j\geq 0$. Hence it follows that
\[
|x_{n_1+j+1}|\leq \| y_{n_0}-x\|_{\infty}+|y_{n_0(n_1+j+1)}|<2\epsilon
\]
for every $j\geq 0$. Thus $x$ belongs to $c_0(\mathbb{N})$. This completes the proof of the theorem.
\end{proof}

\begin{remark}
As we promised in the introduction, we provide below an example of a hypercyclic backward unilateral weighted
shift on the space of square summable sequences $l^{2}(\mathbb{N})$, which is not a $J$-class operator on
$l^{\infty}(\mathbb{N})$. Consider the backward unilateral weighted shift $T$ with weight sequence
\[
(\alpha_1,\alpha_2,\ldots )=
(\frac{1}{2},2,2,\frac{1}{2},\frac{1}{2},2,2,2,\frac{1}{2},\frac{1}{2},\frac{1}{2},2,2,
2,2,\frac{1}{2},\frac{1}{2},\frac{1}{2},\frac{1}{2},\ldots ).
\]
It is easy to check that $T$ is hypercyclic on $l^{2}(\mathbb{N})$. On the other hand we have that
\[
\inf_{j\geq 0}\prod_{i=1}^{n}\alpha_{i+j}\leq\frac{1}{2^n}\quad\textrm{for every}\quad n=1,2,\ldots.
\]
Hence,
\[
\lim_{n\to +\infty} \left( \inf_{j\geq 0}\prod_{i=1}^{n}\alpha_{i+j}\right) =0.
\]
Theorem \ref{t1} implies that $T$ is not a $J$-class operator on $l^{\infty}(\mathbb{N})$.

To complete our study on $J$-class backward unilateral weighted shifts we would like to mention the following
result from \cite{cm}: \textit{a backward unilateral weighted shift $T$ is a $J$-class operator on
$l^{p}(\mathbb{N})$ if and only if $T$ is hypercyclic on $l^{p}(\mathbb{N})$, for $1\leq p<+\infty$}. A similar
result holds for bilateral shifts, see \cite{cm}.
\end{remark}

\begin{theorem}\label{t2}
Let $T:l^{\infty}(\mathbb{Z})\to l^{\infty}(\mathbb{Z})$ be a backward bilateral weighted shift with positive weights $(\alpha_n )_{n\in\mathbb{Z}}$. Then $T$ is not
a $J$-class operator.
\end{theorem}
\begin{proof}
Following a similar line of reasoning as in the proof that (i) implies (ii) in Theorem \ref{t1} and using the
vectors $(\ldots,l,l,l,\ldots )\in l^{\infty}(\mathbb{Z})$ for $l=1,2, \ldots $ we conclude that the sequence
$(\alpha_n )_{n\in\mathbb{Z}}$ is bounded from below by a positive number and
\[
\lim_{n\to +\infty} \left( \inf_{j\in\mathbb{Z}}\prod_{i=1}^{n} \alpha_{i+j}  \right) =+\infty \quad\mbox{and}\quad \lim_{n\to +\infty} \left(
\inf_{j\in\mathbb{Z}}\prod_{i=1}^{n} \alpha_{j-i}  \right) =+\infty.
\]
Assume that there exists a non-zero vector $x=(x_j)_{j\in\mathbb{Z}}\in l^{\infty}(\mathbb{Z})$ such that
$J(x)=l^{\infty}(\mathbb{Z})$. Since $x\neq 0$ there is some $j\in\mathbb{Z}$ such that $x_j\neq 0$. By our
assumption $0\in J(x)$ hence there exist a sequence of positive integers $k_n$ and vectors
$y_n=(y_{nm})_{m\in\mathbb{Z}}$ such that
\[
\| y_n -x\|_{\infty} \to 0 \quad \textrm{and}\quad \| T^{k_n}y_n\|_{\infty} \to 0.
\]
Therefore, taking the $-k_n+1+j$-th coordinate of the vector $T^{k_n}y_n$ we conclude that
\[
\left | \left ( \prod_{i=1}^{k_n} \alpha_{j-i}\right )y_{nj} \right | \to 0.
\]
Since $\prod_{i=1}^{k_n} \alpha_{j-i}\to +\infty$ then $\displaystyle{x_j=\lim_{n\to +\infty}y_{nj}=0}$, a contradiction.
\end{proof}

\begin{corollary}\label{c1}
Let $T$ be a backward unilateral (bilateral) weighted shift with weight sequence $(\alpha_n )_{n\in\mathbb{N}}$
($(\alpha_n )_{n\in\mathbb{Z}}$ respectively). The following are equivalent
\begin{enumerate}
\item[(i)] $J(0)=l^{\infty}(\mathbb{N})$ $(J(0)=l^{\infty}(\mathbb{Z}))$.

\item[(ii)] $\displaystyle{\lim_{n\to +\infty} \left( \inf_{j\geq 0}\prod_{i=1}^{n} \alpha_{i+j}  \right)
=+\infty}$, $(\displaystyle{\lim_{n\to +\infty} \left( \inf_{j\in\mathbb{Z}}\prod_{i=1}^{n} \alpha_{i+j} \right) =+\infty})$.
\end{enumerate}
\end{corollary}

\begin{remark}
By the previous corollary and Theorem \ref{t1} it follows that if $T$ is a backward unilateral weighted shift and $J(0)=l^{\infty}(\mathbb{N})$ then $T$ is $J$-class
operator. However, for backward bilateral weighted shifts this is no longer true. For example consider the backward bilateral weighted shift
$T:l^{\infty}(\mathbb{Z})\to l^{\infty}(\mathbb{Z})$ with weight sequence $(\alpha_n )_{n\in\mathbb{Z}}$, $\alpha_n=2$ for $n\geq 1$ and $\alpha_n=1$ for $n\leq 0$.
Corollary \ref{c1} gives that $J(0)=l^{\infty}(\mathbb{Z})$ and Theorem \ref{t2} implies that $T$ is not a $J$-class operator.
\end{remark}

Using similar arguments as in the proof of Theorem \ref{t1} we obtain the following.

\begin{theorem}\label{t3}
Let $T:l^{\infty}(\mathbb{N})\to l^{\infty}(\mathbb{N})$ be a backward unilateral weighted shift with positive
weights $(\alpha_n )_{n\in\mathbb{N}}$. The following are equivalent.
\begin{enumerate}
\item[(i)] $T$ is a $J^{mix}$-class operator.

\item[(ii)] $\displaystyle{\lim_{n\to +\infty} \left( \inf_{j\geq 0}\prod_{i=1}^{n} \alpha_{i+j}  \right)
=+\infty.}$
\end{enumerate}
In addition, if $T$ is a $J^{mix}$-class operator we have the following complete description of the set of
$J^{mix}$-vectors.
\[
\{ x\in l^{\infty}(\mathbb{N}): J^{mix}(x)=l^{\infty}(\mathbb{N}) \}=c_0(\mathbb{N}).
\]
\end{theorem}

Combining Theorems \ref{t1} and \ref{t3} we obtain the following.

\begin{corollary}
Let $T:l^{\infty}(\mathbb{N})\to l^{\infty}(\mathbb{N})$ be a backward unilateral weighted shift with positive weights $(\alpha_n )_{n\in\mathbb{N}}$. The following
are equivalent.
\begin{enumerate}
\item[(i)] $T$ is a $J^{mix}$-class operator.

\item[(ii)] $T$ is a $J$-class operator.
\end{enumerate}
\end{corollary}

\subsection*{Acknowledgment} We are indebted to the referee for an extremely careful reading of the
manuscript and for sending us a list of very helpful comments concerning the presentation of the paper. He/she
pointed out a gap in the proof of Theorem \ref{t1} and he/she also corrected the statement of Theorem \ref{t2}.

\end{document}